# The Formalization of The Arithmetic System on The Ground of The Atomic Logic

T. J. Stepień, L. T. Stepień

*The Pedagogical University of Cracow, ul. Podchorazych 2, 30 - 084 Krakow, Poland.*

**Abstract:** This paper is a continuation of the paper [16]. Namely, in [16] we have introduced, among others, the definition of the atomic entailment and we have constructed the system $\overset{\sqcap}{S}$, which is based on the atomic entailment. In this paper we show that the classical Arithmetic can be based on the Atomic Logic (see [17]).

## 1. Notational Preliminaries

We assume that $\rightarrow, \sim, \vee, \wedge, \equiv$ denote the connectives of implication, negation, disjunction, conjunction and equivalence, respectively. We use $\Rightarrow, \neg, \Leftrightarrow, \&, \mathbb{V}, \forall, \exists$ as metalogical symbols. Next,

$At_0 = \{p, p_1, p_2, \ldots, q, q_1, q_2, \ldots, s, s_1, s_2, \ldots t, \ldots\}$ denotes the set of all propositional variables. $S_0$ is the set of all well-formed formulas, which are built in the usual manner from propositional variables and by means of logical connectives. $P_0(\phi)$ denotes the set of all propositional variables occuring in $\phi (\phi \in S_0)$. $R_{S_0}$ denotes the set of all rules over $S_0$. Hence, for every $r \in R_{S_0}, \langle \Pi, \phi \rangle \in r$, where $\Pi \subseteq S_0$ and $\phi \in S_0$ and $\Pi$ is a set of premisses and $\phi$ is a conclusion. Hence, $r_*^0$ denotes here the rule of simultaneous substitution for propositional variables. $\langle \{\phi\}, \psi \rangle \in r_*^0 \Leftrightarrow [h^e(\phi) = \psi]$, where $h^e$ is the extension of the mapping $e: At_0 \rightarrow S_0$ $(e \in \varepsilon_*^0)$ to endomorphism $h^e: S_0 \rightarrow S_0$, where

1. $h^e(\phi) = e(\phi)$, for $\phi \in At_0$
2. $h^e(\sim\phi) = \sim h^e(\phi)$
3. $h^e(\phi F \psi) = h^e(\phi) F h^e(\psi)$, for $F \in \{\rightarrow, \vee, \wedge, \equiv\}$ and for every $\phi, \psi \in S_0$.

**Corresponding author:** L. T. Stepień, The Pedagogical University of Cracow, Kraków, Poland. E-mail: sfstepie@cyf-kr.edu.pl, http://www.ltstepien.up.krakow.pl.

Thus, $\varepsilon_*^0$ is a class of functions $e: At_0 \rightarrow S_0$ (for details, see [5], cf. [2]). $r_0^0$ denotes here the Modus Ponens rule in propositional calculus. $R_{0*} = \{r_0^0, r_*^0\}$ (for details, see [2], [5]). A logical matrix is a pair $\mathfrak{M} = \{U, U'\}$, $U$ is an abstract algebra and $U'$ is a subset of the universe $U$, i.e. $U' \subseteq U$. Any $a \in U'$ is called a distinguished element of the matrix $\mathfrak{M}$. $E(\mathfrak{M})$ is the set of all formulas valid in the matrix $\mathfrak{M}$. $\mathfrak{M}_2$ denotes the classical two-valued matrix.

Hence, $Z_2$ is the set of all formulas valid in the classical matrix $\mathfrak{M}_2$ (see [2], [5]).

The symbols $x_1, x_2, \ldots$ are individual variables. $a_1, a_2, \ldots$ are individual constants. $V$ is the set of all individual variables. $C$ is the set of all individual constants. $P_i^n (i, n \in \mathcal{N} = \{1, 2, \ldots\})$ are $n$-ary predicate letters. The symbols $f_i^n (i, n \in \mathcal{N})$ are $n$-ary function letters. The symbols $\bigwedge x_k, \bigvee x_k$ are quantifiers. $\bigwedge x_k$ is the universal quantifier and $\bigvee x_k$ is the existential quantifier. The function letters, applied to the individual variables and individual constants, generate terms. The symbols $t_1, t_2, \ldots$ are terms. $T$ is the set of all terms. $V \cup C \subseteq T$.

The predicate letters, applied to terms, yield simple formulas, i.e. if $P_i^k$ is a predicate letter and $t_1, \ldots, t_k$ are terms, then $P_i^k(t_1, \ldots, t_k)$ is a simple formula.

$Smp$ is the set of all simple formulas. Next, $At_1$ is the set of all atomic formulas, where $At_1 = \{P_i^k(x_{j_1}, \ldots, x_{j_k}): k, i, j_1, \ldots, j_k \in \mathcal{N}\}$. At last, $S_1$ is the set of all well-formed formulas. $FV(\phi)$ denotes the set of all free variables occuring in $\phi$, where $\phi \in S_1$. $x_k \in Ff(t_m, \phi)$ expresses that $x_k$ is free for term $t_m$ in $\phi$. By $x_k/t_m$ we denote the substitution of the term $t_m$ for the individual variable $x_k$. $P_1(\phi)$ denotes the set of all predicate letters occuring in $\phi$ ($\phi \in S_1$). If $FV(\phi) = \{x_1, \ldots, x_k\}$, then $\wedge \phi = \wedge x_1 \ldots \wedge x_k \phi$.

$R_{S_1}$ denotes the set of all rules over $S_1$. Hence, for every $r \in R_{S_1}$, $\langle \Pi, \phi \rangle \in r$, where $\Pi \subseteq S_1$ and $\phi \in S_1$ and $\Pi$ is a set of premisses and $\phi$ is a conclusion. Hence, $r_*^1$ denotes here the rule of simultaneous substitution for predicate letters. $\langle \{\phi\}, \psi \rangle \in r_*^1 \Leftrightarrow [h^e(\phi) = \psi]$, where $h^e$ is the extension of the mapping $e: Smp \to S_1$ ($e \in \varepsilon_*^1$) to endomorphism $h^e: S_1 \to S_1$, where

1. $h^e(\phi) = e(\phi)$, for $\phi \in Smp$
2. $h^e(\sim\phi) = \sim h^e(\phi)$
3. $h^e(\phi F \psi) = h^e(\phi) F h^e(\psi)$, for $F \in \{\to, \vee, \wedge, \equiv\}$
4. $h^e(\wedge x_k \phi) = \wedge x_k h^e(\phi)$
5. $h^e(\vee x_k \phi) = \vee x_k h^e(\phi)$ for every

$\phi, \psi \in S_1$ and $k \in \mathcal{N}$ (for details, see [6], [7]).

Next, $r_0^1$ denotes the Modus Ponens rule in predicate calculus, $r_+^1$ denotes the generalization rule. $R_{0+} = \{r_0^1, r_+^1\}$, $R_{0*+} = \{r_0^1, r_*^1, r_+^1\}$. We write $X \subset Y$, if $X \subseteq Y$ and $X \neq Y$.

We assume here that for every $\alpha \in S_1$, if $FV(\alpha) = \{x_1, \ldots, x_n\}$, then $\alpha^* = \vee x_1 \ldots \vee x_n \sim \alpha$. Hence, for every $\alpha \in S_1$, if $FV(\alpha) = \emptyset$, then $\alpha^* = \sim \alpha$. Analogically, for every $\alpha \in S_0$, $\alpha^* = \sim \alpha$.

Finally, for any $X \subseteq S_i$ and $R \subseteq R_{S_i}$, $Cn_i(R, X)$ is the smallest subset of $S_i$, containing $X$ and closed under the rules $R \subseteq R_{S_i}$, where $i \in \{0,1\}$. The couple $\langle R, X \rangle$ is called a system, whenever $R \subseteq R_{S_i}$ and $X \subseteq S_i$ and $i \in \{0,1\}$. $Syst \cap A_0$ denotes here the class of all systems $\langle R, X \rangle$, which are based on an atomic entailment, where $R \subseteq R_{S_0}$ and $X \subseteq S_0$.

$Syst \cap A_1$ denotes here the class of all systems $\langle R, X \rangle$, which are based on an atomic entailment, where $R \subseteq R_{S_1}$ and $X \subseteq S_1$. $Syst \cap C_1$ denotes here the class of all systems $\langle R, X \rangle$, which are based on a classical entailment, where $R \subseteq R_{S_1}$ and $X \subseteq S_1$.

$\phi \Big|_{R,X}^{A_0} \psi$ denotes that $\psi$ results atomically from $\phi$, on the ground of the system $\langle R, X \rangle$, where $R \subseteq R_{S_0}$ and $X \subseteq S_0$. Next, $\phi \Big|_{R,X}^{A_1} \psi$ denotes that $\psi$ results atomically from $\phi$, on the ground of the system $\langle R, X \rangle$, where $R \subseteq R_{S_1}$ and $X \subseteq S_1$. At last, $\phi \Big|_{R,X}^{C_1} \psi$ denotes that $\psi$ results classically from $\phi$, on the ground of the system $\langle R, X \rangle$, where $R \subseteq R_{S_1}$ and $X \subseteq S_1$ (see [16]).

**Definition 1.1.** *The function* $j: S_1 \to S_0$, *is defined, as follows* (see [16]):
(1) $j(P_k^n(t_1, \ldots, t_n)) = p_k (p_k \in At_0)$
(2) $j(\sim\phi) = \sim j(\phi)$
(3) $j(\phi F \psi) = j(\phi) F j(\psi)$, for $F \in \{\to, \vee, \wedge, \equiv\}$
(4) $j(\wedge x_k \phi) = j(\vee x_n \phi) = j(\phi)$.

## 2. Classical Entailment

**Definition 2.1.** *Let* $Cn_1(R, X) = L \neq \emptyset$ *and* $\phi, \psi \in S_1$. *Then* $\phi \Big|_{R,X}^{C_1} \psi$ *iff the following conditions are satisfied*, [14], [16]:
(1) $(\forall e \in \varepsilon_*^1)[h^e(\wedge \phi) \in L \Rightarrow h^e(\psi) \in L]$
(2) $(\forall e \in \varepsilon_*^1)[h^e((\psi^* \to \phi^*) \to \phi^*) \in L \Rightarrow h^e(\phi^*) \in L]$.

**Definition 2.2.** $\langle R, X \rangle \in Syst \cap C_1$ *iff the following condition is satisfied*, [14], [16]:

$$(\forall \phi, \psi \in S_1)\left[\wedge \phi \to \psi \in Cn_1(R, X) \Leftrightarrow \phi \Big|_{R,X}^{C_1} \psi\right].$$

## 3. The Classical Logic

Let $L_2$ denote the set of all formulas valid in the classical predicate calculus.

Thus, (cf. [6] pp. 68 – 74):

**Theorem 3.1.** $Cn_1(R_{0*+}, L_2) = L_2$.

## 4. Atomic Entailment

In [14], [15] and [16], we have introduced the following definitions (cf. [12], [13]):

**Definition 4.1.** Let $Cn_0(R,X) = L \neq \emptyset$ and $\phi, \psi \in S_0$. Then $\phi \left|\frac{A_0}{R,X}\right. \psi$ iff the following conditions are satisfied:

(1) $(\forall e \in \varepsilon_*^0)[h^e(\phi) \in L \Rightarrow h^e(\psi) \in L$
  $\& P_0(h^e(\phi)) \subseteq P_0(h^e(\psi))]$
(2) $(\forall e \in \varepsilon_*^0)[h^e((\psi^* \to \phi^*) \to \phi^*) \in L \Rightarrow$
  $h^e(\phi^*) \in L \& P_0(h^e(\psi^*)) \subseteq P_0(h^e(\phi^*))]$.

**Definition 4.2.** $\langle R,X \rangle \in Syst \cap A_0$ iff the following condition is satisfied:

$(\forall \phi, \psi \in S_0)\left[\phi \to \psi \in Cn_0(R,X) \Leftrightarrow \phi \left|\frac{A_0}{R,X}\right. \psi\right]$.

**Definition 4.3.** Let $Cn_1(R,X) = L \neq \emptyset$ and $\phi, \psi \in S_1$. Then $\phi \left|\frac{A_1}{R,X}\right. \psi$ iff the following conditions are satisfied:

(1) $(\forall e \in \varepsilon_*^1)[h^e(\wedge\phi) \in L \Rightarrow h^e(\psi) \in L$
  $\& P_1(h^e(\wedge\phi)) \subseteq P_1(h^e(\psi))]$
(2) $(\forall e \in \varepsilon_*^1)[h^e((\psi^* \to \phi^*) \to \phi^*) \in L \Rightarrow$
  $h^e(\phi^*) \in L \& P_1(h^e(\psi^*)) \subseteq P_1(h^e(\phi^*))]$.

**Definition 4.4.** $\langle R,X \rangle \in Syst \cap A_1$ iff the following condition is satisfied:

$(\forall \phi, \psi \in S_1)\left[\wedge\phi \to \psi \in Cn_1(R,X) \Leftrightarrow \phi \left|\frac{A_1}{R,X}\right. \psi\right]$.

## 5. The Atomic Logic

Let us take the matrix (see [16])
$\mathfrak{M}_D = \langle \{0,1,2\}, \{1,2\}, f_D^\to, f_D^\equiv, f_D^\vee, f_D^\wedge, f_D^\sim \rangle$, where:

| $f_D^\to$ | 0 | 1 | 2 |
|---|---|---|---|
| 0 | 1 | 1 | 1 |
| 1 | 0 | 1 | 0 |
| 2 | 0 | 1 | 2 |

| $f_D^\equiv$ | 0 | 1 | 2 |
|---|---|---|---|
| 0 | 1 | 0 | 0 |
| 1 | 0 | 1 | 0 |
| 2 | 0 | 0 | 2 |

| $f_D^\vee$ | 0 | 1 | 2 |
|---|---|---|---|
| 0 | 0 | 1 | 0 |
| 1 | 1 | 1 | 1 |
| 2 | 0 | 1 | 2 |

| $f_D^\wedge$ | 0 | 1 | 2 |
|---|---|---|---|
| 0 | 0 | 0 | 0 |
| 1 | 0 | 1 | 1 |
| 2 | 0 | 1 | 2 |

| $f_D^\sim$ | |
|---|---|
| 0 | 1 |
| 1 | 0 |
| 2 | 2 |

It should be noticed here that the matrix $\mathfrak{M}'_D = \langle \{0,1,2\}, \{1,2\}, f_D^\to, f_D^\sim \rangle$ was investigated by B. Sobocinski (see [9], [10]).

Next, we define $T_D$, putting:

**Definition 5.1.** $T_D = E(\mathfrak{M}_D)$.

In [10] (see [11]) we have proved the following:

**Theorem 5.2.** *The system* $\langle R_{0*}, T_D \rangle$ *is axiomatizable.*

**Theorem 5.3.** *Let* $\phi, \psi \in S_0$ *and*
$(\exists e \in \varepsilon_*^0)[h^e(\phi) \in T_D]$.
*Then* $\phi \to \psi \in Cn_0(R_{0*}, T_D)$ *iff*
  $(\forall e \in \varepsilon_*^0)[h^e(\phi) \in T_D \Rightarrow h^e(\psi) \in T_D$
    $\& P_0(h^e(\phi)) \subseteq P_0(h^e(\psi))]$.

Now we define the system $\hat{S}$ of the Atomic Logic, as follows:

**Definition 5.4.** $\hat{S} = \langle R_0, T_D \rangle$.

Next, we define the set $L_D$, putting (see [16]):

**Definition 5.5.** $L_D = \{\phi \in L_2 : j(\phi) \in T_D\}$.

The system $\overset{\sqcap}{S}$ of the Atomic Logic, is defined, as follows (see [16]):

**Definition 5.6.** $\overset{\sqcap}{S} = \langle R_{0+}, L_D \rangle$.

## 6. The Fundamental Properties of the Systems $\langle R_0, T_D \rangle$, $\langle R_{0+}, L_D \rangle$ and $\langle R_{0+}, L_2 \rangle$

In [16] we have proved the following:

**Theorem 6.1.** $\langle R_0, T_D \rangle \in Syst \cap A_0$.
**Theorem 6.2.** $\langle R_{0+}, L_D \rangle \in Syst \cap A_1$.
**Theorem 6.3.** $\langle R_{0+}, L_2 \rangle \in Syst \cap C_1$.

## 7. The Main Result

**Arithmetic terminology.**

Let $S_A$ denote the set of all well-formed formulas of the Arithmetic System. $FV_A(\phi)$ denotes the set of all free variables occuring in $\phi$, where $\phi \in S_A$. Hence, $\overline{S}_A = \{\phi \in S_A : FV_A(\phi) = \emptyset\}$. $Pr(\phi)$ denotes the set of all predicate letters occuring in $\phi$, where $\phi \in S_A$. $R_{S_A}$ denotes the set of all rules over $S_A$. For any $X \subseteq S_A$ and for any $R \subseteq R_{S_A}$, $Cn(R,X)$ is the smallest subset of $S_A$, containing $X$ and closed under the rules of $R$. The couple $\langle R, X \rangle$ is called a system, whenever $R \subseteq R_{S_A}$ and $X \subseteq S_A$. $R_{0+}^P = \{r_0^P, r_+^P\}$, where $\{r_0^P, r_+^P\} \subseteq R_{S_A}$. $r_0^P, r_+^P$ are Modus Ponens and generalization rule in the Arithmetic System, respectively. Next (cf. [3], [4], [8]):

(1) $\psi^1 : \wedge x_1\, x_1 = x_1$
(2) $\psi^2 : \wedge x_1 \wedge x_2 (x_1 = x_2 \rightarrow x_2 = x_1)$
(3) $\psi^3 : \wedge x_1 \wedge x_2 \wedge x_3 (x_1 = x_2 \rightarrow (x_2 = x_3 \rightarrow x_1 = x_3))$
(4) $\psi^4 : \wedge x_1 \wedge x_2 \wedge x_3 \wedge x_4 (x_1 = x_2 \rightarrow (x_3 = x_4 \rightarrow (x_1 + x_3 = x_2 + x_4)))$
(5) $\psi^5 : \wedge x_1 \wedge x_2 \wedge x_3 \wedge x_4 (x_1 = x_2 \rightarrow (x_3 = x_4 \rightarrow (x_1 \cdot x_3 = x_2 \cdot x_4)))$
(6) $\psi^6 : \wedge x_1 \wedge x_2 \wedge x_3 \wedge x_4 (x_1 = x_2 \rightarrow (x_3 = x_4 \rightarrow (x_1 < x_3 \rightarrow x_2 < x_4)))$
(7) $\psi^7 : \wedge x_1 \sim (1 = x_1 + 1)$
(8) $\psi^8 : \wedge x_1 \wedge x_2 (x_1 + 1 = x_2 + 1 \rightarrow x_1 = x_2)$
(9) $\psi^9 : \wedge x_1 \wedge x_2 (x_1 + (x_2 + 1) = (x_1 + x_2) + 1)$
(10) $\psi^{10} : \wedge x_1 (x_1 \cdot 1 = x_1)$
(11) $\psi^{11} : \wedge x_1 \wedge x_2 [x_1 \cdot (x_2 + 1) = (x_1 \cdot x_2) + x_1]$
(12) $\psi^{12} : \wedge x_1 \wedge x_2 [x_1 < x_2 \equiv \vee x_3\, (x_1 + x_3 = x_2)]$
(13) $\psi^{13} : \bigl(A(1) \wedge \wedge x_1(A(x_1) \rightarrow A(x_1 + 1))\bigr) \rightarrow \wedge x_1 A(x_1)$,

where $A(1), A(x_1), A(x_1 + 1) \in S_A$.
Next, $X_P = \{\psi^1, \psi^2, \psi^3, \psi^4, \psi^5, \psi^6, \psi^7, \psi^8, \psi^9, \psi^{10}, \psi^{11}, \psi^{12}\}$.

$Y_P$ denotes here the set of all axioms of induction. At last, $L_2^r$ and $X_r$ denote the set of all logical axioms of the Arithmetic System and the set of all specific axioms of the Arithmetic System, where $L_2^r, X_r \subseteq S_A$, respectively. $Sx$ denotes here the successor of $x$ (see [1]).

$\psi^{14}$ denotes the formula (see [3])
$\wedge x_1 \wedge x_2 [\vee x_3 (Sx_3 + x_1 = x_2) \equiv (x_1 < x_2)]$.

**Definition 7.1.** *The function* $i : S_A \rightarrow S_0$, *is defined, as follows*:
(1) $i(t_n = t_m) = p_k\, (p_k \in At_0)$
(2) $i(t_n < t_m) = p_s\, (p_s \in At_0)$
(3) $i(\sim\phi) = \sim i(\phi)$
(4) $i(\phi F \psi) = i(\phi) F i(\psi)$, for $F \in \{\rightarrow, \vee, \wedge, \equiv\}$
(5) $i(\wedge x_k\, \phi) = i(\vee x_n\, \phi) = i(\phi)$,

where $\phi, \psi \in S_A$.

**Definition 7.2.** $\langle R_{0+}^P, L_2^r \cup X_r \rangle$ *is the Arithmetic System*, where $X_r = X_P \cup Y_P$ (see [3], [4], [8]).

In [8], one can read that the System $\langle R_{0+}^P, L_2^r \cup X_r \rangle$ is a modification of Peano's Arithmetic System.

Next,
**Definition 7.3.** $L_D^r = \{\phi \in L_2^r : i(\phi) \in T_D\}$.

**Theorem 7.4.**
$Cn(R_{0+}^P, L_D^r \cup X_r) = Cn(R_{0+}^P, L_2^r \cup X_r)$,
where $X_r = X_P \cup Y_P$.

Proof. Let
(1) $\alpha \in L_2^r - L_D^r$
and
(2) $X_r = X_P \cup Y_P$.

From **Definition 5.1.**, **Definition 5.5.**, **Definition 7.1.** and **Definition 7.3.**, it follows that
(3) $(\forall \phi \in L_2^r)\, [Pr(\phi) \subseteq \{=\} \Rightarrow \phi \in L_D^r]$
and
(4) $(\forall \phi \in L_2^r)\, [Pr(\phi) \subseteq \{<\} \Rightarrow \phi \in L_D^r]$.

From (1), (3) and (4), it follows that
(5) $< \in Pr(\alpha)$
and
(6) $= \in Pr(\alpha)$.

From (1), (3), (4), (5), (6), **Definition 5.1.**, **Definition 5.5.**, **Definition 7.1.**, and the definition of the formula $\psi^{12}$, it follows that

(7) $i(\psi^{12} \to \alpha) \in T_D$

and

(8) $\psi^{12} \to \alpha \in L_2^r$.

Hence, from (1), (5), (6), by **Definition 5.1.**, **Definition 5.5.** and **Definition 7.3.**, it follows that

(9) $\psi^{12} \to \alpha \in L_D^r$.

Hence, from (1), (2), (3), (4), **Definition 7.2.** and **Definition 7.3.**, it follows that

(10) $Cn(R_{0+}^P, L_D^r \cup X_r) = Cn(R_{0+}^P, L_2^r \cup X_r)$.  □

Thus, by the proof of **Theorem 7.4.**, one can obtain (see [17]):

**Conclusion 7.5.**

*Every Arithmetic System* $\langle R_{0+}^P, L_2^r \cup X_r \rangle$ *can be based on the system of the Atomic Logic* $\langle R_{0+}, L_D^r \rangle$, *where*

$$\psi^{12} \in Cn(R_{0+}^P, L_2^r \cup X_r)$$

*or*

$$\psi^{14} \in Cn(R_{0+}^P, L_2^r \cup X_r)$$

*and* $(\forall \phi \in L_2^r \cup X_r)[Pr(\phi) \subseteq \{=, <\}]$,

(cf. [1] p.530 – 541).

## References

[1] B. Buldt. The Scope of Gödel's First Incompleteness Theorem, *Logica Universalis*, 8:499 – 552, 2014.

[2] Yu. L. Ershov and E. A. Palyutin. *Mathematical Logic*. Mir Publishers, Moscow, 1984.

[3] A. Grzegorczyk. *An Outline of Mathematical Logic. Fundamental Results and Notions Explained with All Details*, D. Reidel Publishing Company, Dordrecht-Holland/Boston-USA, PWN, Warszawa, 1974.

[4] R. Murawski. *Recursive Functions And Metamathematics. Problems of Completeness and Decidability, Gödel's Theorems*. Springer Science+Business Media Dordrecht, 1999.

[5] W. A. Pogorzelski. *The Classical Propositional Calculus*. PWN, Warszawa, 1975.

[6] W. A. Pogorzelski. *The Classical Calculus of Quantifiers*. PWN, Warszawa, 1981.

[7] W. A. Pogorzelski and T. Prucnal. The substitution rule for predicate letters in the first-order predicate calculus. *Reports on Mathematical Logic*, 5:77 – 90, 1975.

[8] H. Rasiowa. *Introduction to Modern Mathematics*. North-Holland Publishing Company, 1973.

[9] B. Sobocinski. Axiomatization of partial system of three-valued calculus of propositions. *The Journal of Computing Systems*, 1:23 – 55, 1952.

[10] T. Stepien. System $\overline{S}$. *Reports on Mathematical Logic*, 15:59 – 65, 1983.

[11] T. Stepien. System $\overline{S}$. *Zentralblatt für Mathematik*, 471, 1983.

[12] T. Stepien. Logic based on atomic entailment. *Bulletin Of The Section Of Logic*, 14:65-71, 1985.

[13] T. Stepien. Logic Based On Atomic Entailment And Paraconsistency. *11th International Congress Of Logic, Methodology And Philosophy Of Science* (August 1999, Krakow, Poland).

[14] T. J. Stepien and L. T. Stepien. Atomic Entailment and Classical Entailment, *The Bulletin of Symbolic Logic*, 17:317 – 318, 2011.

[15] T. J. Stepien and L. T. Stepien. Atomic Entailment and Atomic Inconsistency, *6th International Conference "Non-classical logics. Theory & Applications"* (4 – 6 September 2013, Lodz, Poland).

[16] T. J. Stepien and L. T. Stepien. Atomic Entailment and Atomic Inconsistency and Classical Entailment. *Journal of Mathematics and System Science*, 5:60 – 71, 2015.

[17] T. J. Stepien and L. T. Stepien. The formalization of the arithmetic system on the ground of atomic logic. (Logic Colloquium 2015, 3 – 8 August 2015, Helsinki, Finland), to appear in *The Bulletin of Symbolic Logic*.